\magnification 1200
\advance\hoffset by 1,5truecm
\advance\hsize by 2,3 truecm
\def\makefootline{\baselineskip=52pt\line{\the\footline}}
\vsize= 23 true cm
\hsize= 14 true cm
\overfullrule=0mm

\font\grandsy=cmsy10 scaled \magstep0
\def\SS{{\grandsy x}}

\headline={\hfill\tenrm\folio\hfil}
\footline={\hfill}\pageno=1

\newcount\coefftaille \newdimen\taille
\newdimen\htstrut \newdimen\wdstrut
\newdimen\ts \newdimen\tss

\def\fixetaille#1{\coefftaille=#1
\htstrut=8.5pt \multiply \htstrut by \coefftaille \divide \htstrut by 1000
 \wdstrut=3.5pt \multiply \wdstrut by \coefftaille \divide \wdstrut by 1000
\taille=10pt  \multiply \taille by \coefftaille \divide \taille by 1000
\ts=\taille \multiply \ts by 70 \divide \ts by 100
 \tss=\taille \multiply \tss by 50 \divide \tss by 100
\font\tenrmp=cmr10 at \taille
\font\sevenrmp=cmr7 at \ts
\font\fivermp=cmr5 at \tss
\font\tenip=cmmi10 at \taille
\font\sevenip=cmmi7 at \ts
\font\fiveip=cmmi5 at \tss
\font\tensyp=cmsy10 at \taille
\font\sevensyp=cmsy7 at \ts
\font\fivesyp=cmsy5 at \tss
\font\tenexp=cmex10 at \taille
\font\tenitp=cmti10 at \taille
\font\tenbfp=cmbx10 at \taille
\font\tenslp=cmsl10 at \taille}

\def\fspeciale{\textfont0=\tenrmp%
\scriptfont0=\sevenrmp%
\scriptscriptfont0=\fivermp%
\textfont1=\tenip%
\scriptfont1=\sevenip%
\scriptscriptfont1=\fiveip%
\textfont2=\tensyp%
\scriptfont2=\sevensyp%
\scriptscriptfont2=\fivesyp%
\textfont3=\tenexp%
\scriptfont3=\tenexp%
\scriptscriptfont3=\tenexp%
\textfont\itfam=\tenitp%
\textfont\bffam=\tenbfp%
\textfont\slfam=\tenbfp%
\def\it{\fam\itfam\tenitp}%
\def\bf{\fam\bffam\tenbfp}%
\def\rm{\fam0\tenrmp}%
\def\sl{\fam\slfam\tenslp}%
\normalbaselineskip=12pt%
\multiply \normalbaselineskip by \coefftaille%
\divide \normalbaselineskip by 1000%
\normalbaselines%
\abovedisplayskip=10pt plus 2pt minus 7pt%
\multiply \abovedisplayskip by \coefftaille%
\divide \abovedisplayskip by 1000%
\belowdisplayskip=7pt plus 3pt minus 4pt%
\multiply \belowdisplayskip by \coefftaille%
\divide \belowdisplayskip by 1000%
\setbox\strutbox=\hbox{\vrule height\htstrut depth\wdstrut width 0pt}%
\rm}

\def\vmid#1{\mid\!#1\!\mid}

\def\fle{\rightarrow}

\null\vskip-1cm

\font\sc=cmcsc10
\font\grandsy=cmsy10 scaled\magstep2

\newdimen\emm 
\def\pmb#1{\emm=0.03em\leavevmode\setbox0=\hbox{#1}
\kern0.901\emm\raise0.434\emm\copy0\kern-\wd0
\kern-0.678\emm\raise0.975\emm\copy0\kern-\wd0
\kern-0.846\emm\raise0.782\emm\copy0\kern-\wd0
\kern-0.377\emm\raise-0.000\emm\copy0\kern-\wd0
\kern0.377\emm\raise-0.782\emm\copy0\kern-\wd0
\kern0.846\emm\raise-0.975\emm\copy0\kern-\wd0
\kern0.678\emm\raise-0.434\emm\copy0\kern-\wd0
\kern\wd0\kern-0.901\emm}

\font\tendb=msbm10
\font\sevendb=msbm7

\newfam\dbfam
\textfont\dbfam=\tendb\scriptfont\dbfam=\sevendb\scriptscriptfont\dbfam=\sevendb
\def\db{\fam\dbfam\tendb}

\def\C{{\db C }}

\def\N{{\db N }}

\def\R{{\db R }}

\def\cos{\mathop{\rm cos}\nolimits}


\font\gothique=eufm10
\def\got#1{{\gothique #1}}

\newdimen\margeg \margeg=0pt
\def\bb#1&#2&#3&#4&#5&{\par{\parindent=0pt
    \advance\margeg by 1.1truecm\leftskip=\margeg
    {\everypar{\leftskip=\margeg}\smallbreak\noindent
    \hbox to 0pt{\hss\bf [#1]~~}{\bf #2 - }#3~; {\it #4.}\par\medskip
    #5 }
\medskip}}

\newdimen\margeg \margeg=0pt
\def\bbaa#1&#2&#3&#4&#5&{\par{\parindent=0pt
    \advance\margeg by 1.1truecm\leftskip=\margeg
    {\everypar{\leftskip=\margeg}\smallbreak\noindent
    \hbox to 0pt{\hss [#1]~~}{\pmb{\sc #2} - }#3~; {\it #4.}\par\medskip
    #5 }
\medskip}}

\newdimen\margeg \margeg=0pt
\def\bba#1&#2&#3&#4&#5&{\par{\parindent=0pt
    \advance\margeg by 1.1truecm\leftskip=\margeg
    {\everypar{\leftskip=\margeg}\smallbreak\noindent
    \hbox to 0pt{\hss [#1]~~}{{\sc #2} - }#3~; {\it #4.}\par\medskip
    #5 }
\medskip}}

\def\messages#1{\immediate\write16{#1}}

\def\findem{\vrule height0pt width4pt depth4pt}

\long\def\demA#1{{\parindent=0pt\messages{debut de preuve}\smallbreak
     \advance\margeg by 2truecm \leftskip=\margeg  plus 0pt
     {\everypar{\leftskip =\margeg  plus 0pt}
              \everydisplay{\displaywidth=\hsize
              \advance\displaywidth  by -1truecm
              \displayindent= 1truecm}
     {\bf Proof } -- \enspace #1
      \hfill\findem}\bigbreak}\messages{fin de preuve}}

\def\resp{\mathop{\rm resp}\nolimits}
\def\resp.{\mathop{\rm resp.}\nolimits}



\null\vskip 1,5cm
\centerline{\bf Holomorphic Cliffordian Functions}
\bigskip
\centerline{\bf by}
\medskip
\centerline{{\sc Guy Laville } and {\sc Ivan Ramadanoff}}

\bigskip\bigskip
{\leftskip=0cm\rightskip=10pt
{\parindent=0cm\narrower\fixetaille{700}{{\fspeciale {\bf Abstract.-} \  
{\sevenrm The aim of this paper is to put the fundations of a new theory of
functions, called holomorphic Cliffordian, which should play an essential role in
the generalization of holomorphic functions to higher dimensions.
\medskip
Let \ $\R_{0,2m+1}$ \ be the Clifford algebra of \ $\R^{2m+1}$ \ with a quadratic
form of negative signature, \ $D = \displaystyle\sum_{j=0}^{2m+1} \ e_j \
\displaystyle{\partial\over \partial x_j}$\ be the usual operator for monogenic
functions and $\Delta$ the ordinary Laplacian. The holomorphic Cliffordian
functions are functions \ $f : \R^{2m+2} \fle \R_{0,2m+1}$, which are solutions of
$D \Delta^m f = 0$.

\medskip
Here, we will study polynomial and singular solutions of this equation, we will
obtain integral representation formulas and deduce the analogous of the Taylor and
Laurent expansions for holomorphic Cliffordian functions.

\medskip
In a following  paper, we will put the fundations of the Cliffordian elliptic
function theory.} 
}
}
\par}

\bigskip\bigskip
\bigskip\bigskip
{\bf 0. Introduction} 
\bigskip
\hskip0,45cm The classical theory of holomorphic functions of one complex variable
has been genera\-li\-zed in two directions. The first is the theory of holomorphic
functions of several complex variables~:  in this case we keep the field $\C$ and
take the system of partial differential operators \ $\partial / \partial
\overline z_i$, \ $i = 1,\ldots , n$.  The second direction is the theory of
monogenic functions~: in this case we take the Clifford algebra and take the
operator \ $D = \displaystyle\sum_{j=0}^m e_i \ \partial / \partial x_i$ \ ($\{
e_i\}$ orthogonal basis).
\medskip
\hskip0,45cm Here we follow a different path : we think that the most important
thing in the theory of one complex variable is the fact that the identity  (i.e.
$z$)  and its powers  (i.e. $z^n$)  are holomorphic.

\bigskip\bigskip
\bigskip\bigskip
{\bf 1. Notations} 
\bigskip
\hskip0,45cm
Let \ $\R_{0,2m+1}$ be the Clifford algebra of the real vector space $V$ of
dimension $2m+1$, provided with a quadratic form of negative signature, $m\in\N$. 
Denote by $S$ the set of the scalars in \ $\R_{0,2m+1}$,  which can be identified
to \ $\R$.  Let \ $\{ e_i\}, i = 1,2,\ldots , 2m+1$ be an orthonormal basis of\ 
$V$ and let $e_0 = 1$.

A point \ $x = (x_0, x_1,\ldots , x_{2m+1})$ \ of \ $\R^{2m+2}$ \ could be also
considered as an element of $S \oplus V$, namely \ $x =
\displaystyle\sum_{i=0}^{2m+1} e_i x_i$.  So, \ $x$, \ being in \ $S\oplus V$,  is
in the Clifford algebra $\R_{0,2m+1}$  and we can act on him by the principal
involution in \ $\R_{0,2m+1}$,  which will coincide with a  kind of
``conjugation"~:
$$x^* = x_0 - \sum_{i=1}^{2m+1} e_i x_i.$$

It is remarkable that
$$xx^* = x^*x = \ \vmid{x}^2,$$

where \ $\vmid{x}$ \ denotes the usual euclidean norm of $x$ in \ $\R^{2m+2}$.

Sometime, if necessary, we will resort to the notation $x = x_0 + \overrightarrow
x$,  where \ $\overrightarrow x$ \ is the vector part of $x$, namely \
$\overrightarrow x = \displaystyle\sum_{i=1}^{2m+1} e_i x_i$.

\vskip 1,5cm
{\bf 2. General definitions}
\bigskip
\hskip0,45cm Let \ $\Omega$ \ be an open set of \ $S\oplus V$.  We will be
interested in functions \ $f : \Omega \fle \R_{0,2m+1}$.  It should be noted that
one might consider only functions \ $f : \Omega \fle S\oplus V$.  The last ones
generate the previous by means of (right) linear combinations.

It is well known that the following operator, named Cauchy  (or Fueter, or
Dirac) operator ([1], [2], [3], [4]) lies on the basis of the theory of
(left) monogenic functions~:
$$D = \sum_{i=0}^{2m+1} \ e_i \ {\partial\over \partial x_i}.\leqno (1)$$

A function \ $f : \Omega \fle \R_{0,2m+1}$ \ is said to be  (left) monogenic in \
$\Omega$ \ if and only if~:
$$Df (x) = 0$$

for each \ $x$ \ on \ $\Omega$.

It is important to note that the operator \ $D$ \ possesses a conjugate operator \
$D^*$~:
$$D^* = {\partial\over \partial x_0} - \sum_{i=1}^{2m+1} \ e_i \ {\partial\over
\partial x_i},\leqno (2)$$
and that \ $DD^* = D^*D = \Delta$,  where \ $\Delta$ \ is the ordinary Laplacian.

Now let us state the following~:

\bigskip\bigskip
{\bf D\'efinition.-} \ {\sl A function \ $f : \Omega \fle \R_{0,2m+1}$ \ is said
to be (left) holomorphic Cliffordian in \ $\Omega$ \ if and only if~:
$$D \Delta^m f(x) = 0$$

for each \ $x$ \ of \ $\Omega$.  Here \ $\Delta^m$ \ means the $m$ times iterated 
Laplacian $\Delta$.
}

\bigskip\bigskip
{\bf Remark.-} \ The set of holomorphic Cliffordian functions is wider than the
set of monogenic functions in the sense that every monogenic function is also a
holomorphic Cliffordian, but the reciproque is false. Indeed, if \ $Df = 0$, \
then \ $D \Delta^m f = \Delta^m Df = 0$ \ because the operator \ $\Delta^m$ \ is a
scalar operator.

The simplest example of a function which is holomorphic Cliffordian, but not
monogenic is the identity,  \ id : $x\longmapsto x$,  for which \ $Dx = - 2m \not=
0$ \ and clearly \ $D\Delta^mx = 0$.

Later, we will be able to prove that all entire powers of $x$ are holomorphic
Cliffordians, while they are not monogenics.

\bigskip\bigskip
{\bf Remark.-} \ $f$ \ is (left)  holomorphic Cliffordian if and only if   \
$\Delta^mf$ \ is (left) monogenic.

\vskip 1,5cm
{\bf 3. Some properties of the holomorphic Cliffordian functions}
\bigskip
\hskip0,45cm (i) \ All the components of the so called scalar, vector, bivector,
$\ldots$, up to the pseudo-scalar parts of a holomorphic Cliffordian function $f$
are polyharmonics of order $m+1$. This is obvious taking into account that, if
$D\Delta^mf = 0,$ then applying $D^\ast,$ one get $\Delta^{m+1}f = 0$ and the
result follows because $\Delta^{m+1}$ is a scalar operator.

\medskip
\hskip0,45cm (ii)\ If $f$ is a polyharmonic function of order $m+1,$ i.e.
$\Delta^{m+1}f=0,$ then the function $D^\ast f$ is holomorphic Cliffordian. Indeed,
$\Delta^{m+1}f = D D^\ast\Delta^mf = D \Delta^m(D^\ast f) = 0.$

This property will play an important role in the next part of this paper because
it is a good machinery for generating holomorphic Cliffordian functions.

\medskip
\hskip0,45cm (iii)\ Let us compute $\Delta(xg),$ where $g : S\oplus V
\fle\R_{0,2m+1}$ is sufficiently smooth. One has~:
\medskip
$$\eqalign{&\Delta(xg) = \displaystyle\sum_{i=0}^{2m+1} {\partial^2\over\partial
x_i^2}(xg) =
\sum_{i=0}^{2m+1}{\partial\over\partial x_i}\Bigl[\Bigl({\partial\over\partial
x_i}(x_0+\overrightarrow x)\Bigr)g + x {\partial g\over\partial x_i}\Bigr]\cr
&= \displaystyle\sum_{i=0}^{2m+1}{\partial\over\partial x_i}\Bigl(e_ig + x
{\partial g\over\partial x_i}\Bigr) = \sum_{i=0}^{2m+1}\Bigl(e_i {\partial
g\over\partial x_i} + {\partial x\over\partial x_i} {\partial g\over\partial x_i}
+ x {\partial^2g\over\partial x_i^2}\Bigr) = 2 Dg + x\Delta g.\cr}$$
\medskip
Thus, we have :
\medskip
\hskip 3cm $2Dg = \Delta(xg) - x\Delta g$

\hskip 3cm $x\Delta g = \Delta(xg) - 2Dg$
\medskip
Now, if we compute : $2D\Delta g = \Delta(x\Delta g) - x\Delta^2g =
\Delta(\Delta(xg) - 2Dg) - x\Delta^2g \break = \Delta^2(xg) - 2D\Delta
g-x\Delta^2g.$ In this way, we get
$$4D\Delta g = \Delta^2(xg) - x\Delta^2g.$$

\medskip
Using a recurrence process, we obtain :
$$2(p+1)D\Delta^pg = \Delta^{p+1}(xg) - x \Delta^{p+1}g,$$

for every $p\in\N.$ Putting in the last formula, $p=m,$ one deduces~:
\medskip
$$2(m+1)D\Delta^mg = \Delta^{m+1}(xg) - x\Delta^{m+1}g.\leqno (3)$$
\medskip
which gives a sufficient condition for $g$ to be holomorphic Cliffordian, namely
$g$ and $xg$ have to be polyharmonics of order $(m+1).$ 

\medskip
But this condition is also necessary. If $g$ is holomorphic Cliffordian,
$D\Delta^mg = 0$ and, using (3) one has~:
\medskip
$$\Delta^{m+1}(xg) = x \Delta^{m+1}g.\leqno (4)$$
\medskip
Now, compute the right hand side : $x\Delta^{m+1}g = x D^\ast(D\Delta^mg) = 0.$
So, $xg$ is polyharmonic. From (4), again, it follows that $g$ is also
polyharmonic.

\bigskip
\hskip0,45cm (iv) \  The equation \ $D\Delta^mf = 0$ \ is equivalent to the system
$$\left\{ \eqalign{
&Df_{(2p+1)} = f_{(2p+2)}\cr
&D^*f_{(2p+2)} = f_{(2p+3)}\cr
&Df_{(2m+1)} = 0. \cr}\right.$$

with \ \ $p = 0, 1,\ldots , m-1$ \  and \  $f_{(1)} = f$.

\vfill\eject
{\bf 4. Some examples of holomorphic Cliffordian functions}
\bigskip
\hskip0,45cm Let us start with the following lemma~:

\bigskip\bigskip
{\bf Lemma.-}\ If $f : \R^2 \fle\R$,  is harmonic, then
$f(x_0,\vmid{\hskip-0,08cm\vec x}),$ where $x = x_0+\vec x$ and\break
$\vmid{\vec x}^2 = \displaystyle\sum_{i=1}^{2m+1}x_i^2$, \  is
$(m+1)-$harmonic, that is~:
$$\Delta^{m+1}f(x_0,\vmid{\vec x}) = 0.$$

\bigskip\bigskip
\demA{
Set $r = \vmid{\vec x}.$ Thus, the Laplacian could be written as :
$$\Delta = {\partial^2\over\partial x_0^2} + {\partial^2\over\partial r^2} +
{2m\over r}{\partial\over\partial r}.$$

But $f(x_0,r)$ is harmonic, so \  $\displaystyle{\partial^2f\over\partial x_0^2} +
{\partial^2f\over\partial r^2} = 0$\  and hence~:
$$\Delta f(x_0,\vmid{\vec x}) = {2m\over r}{\partial f\over\partial
r}.$$

Now, compute the first iteration :
$${1\over 2m} \Delta^2 f(x_0,\vmid{\vec x}) = {1\over r}
{\partial^3f\over\partial x_0^2\partial r} + {2\over r^3} {\partial f\over\partial
r} - {2\over r^2} {\partial^2f\over\partial r^2}$$

$$ + {1\over r}{\partial^3f\over\partial r^3} + {2m\over r}\Bigl(-{1\over
r^2}{\partial f\over\partial r} + {1\over r}{\partial^2f\over\partial r^2}\Bigr)
=$$

$$= {2m-2\over r^2} {\partial^2f\over\partial r^2} - {2m-2\over r^3} {\partial
f\over\partial r}.$$

Here, we have take into account that $\displaystyle{1\over
r}{\partial\over\partial r}\Bigl({\partial^2f\over\partial x_0^2} +
{\partial^2f\over\partial r^2}\Bigr) = 0.$ Thus, we get~:
$${1\over 2m} \cdot {1\over 2m-2} \Delta^2 f(x_0,\vmid{\vec x}) =
{1\over r^2} {\partial^2f\over\partial r^2} - {1\over r^3}{\partial f\over\partial
r}$$

It is easy to show that :
$${1\over 2m}\cdot{1\over 2m-2} \Delta^2 f(x_0,\vmid{\vec x}) =
\Bigl({1\over r}{\partial\over\partial r}\Bigr)^2 f.$$

Using a recurrence process, it is possible to prove that
$$\Delta^k f(x_0, \ \vmid{\vec x}) = 2m (2m-2) \cdots (2m-2k+2) \ {\Bigl(
{1\over r} \ {\partial\over \partial r}\Bigr)}^k f, \leqno\hskip 1cm (5)$$
for \ $k\in\N$.  In fact, one needs also a preliminary formula~:
$${\partial^2\over \partial r^2} \ {\Bigl( {1\over r} \ {\partial\over \partial
r}\Bigr)}^k = - 2k \ {\Bigl( {1\over r} \ {\partial\over \partial r}\Bigr)}^{k+1}
+ {\Bigl( {1\over r} \ {\partial\over \partial r}\Bigr)}^k \ {\partial^2\over
\partial r^2}$$
the proof of which is also achieved by a reccurence argument. The end of the proof
of the lemma would be performed setting in  (5), $k = m+1$.
}

\bigskip
\hskip0,45cm Now, combining this lemma with the property (ii), we get a
nice process for generating holomorphic Cliffordian functions. Let us illustrated
this by the following~:

\bigskip\bigskip
{\bf Proposition.-} \ Let \ $x = x_0 + \vec x = x_0 +
\displaystyle\sum_{i=1}^{2m+1} \ e_i x_i, \lambda\in\R$ and $n\in\N$. Then, the
functions \ $x\longmapsto e^{\lambda x}$ \ and \ $x \longmapsto x^n$ \ are
holomorphic Cliffordians.

\bigskip
\demA{It is clear that it suffices to prove that \ $\Delta^m e^{\lambda
x}$ \ and \ $\Delta^m x^n$ \ are monogenics.  By the lemma, taking the real part of
$e^{\lambda z}$,  where $z\in\C$,  one has~:
$$\Delta^{m+1} \ e^{\lambda x_{0}} \ \cos (\lambda \ \vmid{\vec x}) \  =
0.$$

We will obtain a holomorphic Cliffordian function taking \ $D^* \ e^{\lambda
x_{0}} \cos (\lambda \vmid{\vec x})$.  Let us compute this~:
\medskip

$D^\ast e^{\lambda x_0}\cos (\lambda\vmid{\vec x})\  = \lambda
e^{\lambda x_0} \cos (\lambda\vmid{\vec x}) - \lambda e^{\lambda
x_0}\sin (\lambda\vmid{\vec x}) D^\ast(\vmid{\vec x})$

$= \lambda e^{\lambda x_0}\Bigl[\cos (\lambda\vmid{\vec
x}) -\sin (\lambda\vmid{\vec x}) \displaystyle{D^\ast(\vmid{\vec
x}^2)\over 2\vmid{\vec x}}\Bigr] =$

$\displaystyle = \lambda e^{\lambda x_0}\Bigl[\cos (\lambda\vmid{\vec x}) +
{\vec x\over\vmid{\vec x}}\sin (\lambda\vmid{\vec
x})\Bigr] = \lambda e^{\lambda x_0}e^{\lambda\vec x} = \lambda
e^{\lambda x}.$

\bigskip
\hskip0,45cm It follows immediatly that all the terms of the expansion of \
$e^{\lambda x}$ \ are holomorphic Cliffordian, and in particular \ $x^n$ \ for
$n\in\N$.
}

\bigskip\bigskip
{\bf Remark.-} When $f$ is holomorphic Cliffordian, then the same is true for
all\break 
$\displaystyle{\partial\over\partial x_j} f$, \quad $j=0,\ldots,2m+1.$ Indeed~:
$$D\Delta^m\Bigl({\partial\over\partial x_j}f) = {\partial\over\partial
x_j}(D\Delta^mf) = 0.$$

\hskip0,45cm More generally, let us denote by $D^\alpha$ the operator of
derivation~:
$$D^\alpha = {\partial^{\alpha_0+\alpha_1+\ldots+\alpha_{2m-1}}\over\partial
x_0^{\alpha_0}\partial x_1^{\alpha_1}\ldots\partial x_{2m+1}^{\alpha_{2m+1}}}.$$

\hskip0,45cm Where $\alpha = (\alpha_0,\alpha_1,\ldots\alpha_{2m+1}) \in\N^{2m+2}$
is a multiindice, then if $f$ is holomorphic Cliffordian, then $D^\alpha f$ est also
holomorphic Cliffordian.

\hskip0,45cm See [2] and [4].

\vskip 1,5cm
{\bf 5. Polynomial solutions of $\bf D \Delta^m f = 0$}
\bigskip
\hskip0,45cm Now, we know that all integer powers of $x$ are
monomials which are solutions of the equation
$$D \Delta^m(x^n) = 0, \qquad n\in\N.$$

\hskip0,45cm Let us find all possible "monomials``. For this purpose, set
$\alpha = (\alpha_0,\alpha_1,\ldots,\alpha_{2m+1})$ with $\alpha_i\in\N$ and
$\vmid\alpha = \displaystyle\sum_{i=0}^{2m+1}\alpha_i.$ Consider the set
$\{e_\nu\} =
\{e_0,\ldots,e_0,e_1,\ldots,e_1,\break\ldots,e_{2m+1},\ldots,e_{2m+1}\}$ where $e_0$
is written $\alpha_0$ times, $e_i : \alpha_i$ times and $e_{2m+1} : \alpha_{2m+1}$
times. Then set~:
$$P_\alpha(x) = \sum_{\hbox{\got S}} \prod_{\nu=1}^{\vmid\alpha-1}
(e_{\sigma(\nu)}x)e_{\sigma(\vmid\alpha)},\leqno (6)$$

the sum being expanded over all distinguishable elements $\sigma$ of the permutation
group
$\hbox{\got S}$ of the set $\{e_\nu\}.$

\hskip0,45cm The function $P_\alpha(x),$ as a function of $x,$ is a polynomial of
degree $\vmid\alpha-1.$ A straightforward calculation carried on
$P_\alpha$ shows that $P_\alpha$ is equal, up to a rational constant, to
$D^\alpha(x^{2\vmid\alpha-1}).$ It follows then that $P_\alpha(x)$ is a
holomorphic Cliffordian function.

\hskip0,45cm As an illustration, let us compute $P_{(0,1,1,0)}(x), P_{(1,1,0,0)}(x)$
and $P_{(2,0,0,0)}(x)$ in the case when $m=1.$ Following our notations, we have \
$\vmid{\alpha} \ = 2$ and
$$\eqalign{
&P_{(0,1,1,0)}(x) = e_1x e_2+ e_2xe_1 \cr
&P_{(1,1,0,0)}(x) = e_0x e_1+ e_1xe_0 \cr
&P_{(2,0,0,0)}(x) = e_0xe_0.\cr}$$

\hskip0,45cm Now, as the first polynomial is concerned, let us calculate
$$\displaylines{
{\partial^2\over \partial x_1 \partial x_2} (x^3) = {\partial\over \partial x_1} \
(e_2x^2 + xe_2x + x^2e_2) = \cr
= (e_2e_1x + e_2xe_1) + (e_1e_2 x + xe_2e_1) + (e_1xe_2 + xe_1e_2) \cr
= e_1xe_2 + e_2xe_1 = P_{(0,1,1,0)}(x) \cr}$$

For the second one :
$${\partial^2\over \partial x_0\partial x_1} (x^3) = {\partial\over \partial x_1}
(3x^2) = 3 (e_1x+xe_1) = 3 (e_0xe_1 + e_1xe_0) = 3 \ P_{(1,1,0,0)} (x).$$

Finally :
$${\partial^2\over \partial x_0^2} (x^3) = 6x = 6 \ e_0xe_0 = 6 \
P_{(2,0,0,0)}(x).$$

The general formula is :
$$D^\alpha (x^{2\vmid{\alpha}-1}) = \cases{
P_\alpha (x) , &\hbox{if} \ $\alpha_0 = 0$\cr
\alpha_0! \ C_{2\vmid{\alpha}-1}^{\alpha_{0}} P_\alpha (x), &\hbox{if} \ $\alpha_0
\not= 0.$\cr}$$

\bigskip\bigskip
Later, we will be able to prove that the polynomials $P_\alpha (x)$ form a basis
of the space of polynomial solution of the equation $D\Delta^mP = 0$.

\bigskip
{\bf Remark :} \ the polynomials \ $P_\alpha (x)$ \ are left and right holomorphic
Cliffordian.

Put
$$\eqalign{
&\lambda = \sum_{i=0}^{2m+1} \ \lambda_i e_i , \quad \lambda_i \in \R \cr
&\lambda_\alpha = \prod_{i=1}^{2m+1} \ \lambda_i^{\alpha_{i}} \cr}$$

then, the following formal series gives the generating function~:
$$(1 - \lambda x)^{-1} \lambda = \sum_\alpha \ P_\alpha (x) \lambda_\alpha .$$

It is convenient, for certain computations, to modify slightly these polynomials 
a little bit~:

Let \quad $\overrightarrow\alpha = (\alpha_1,\ldots , \alpha_{2m+1})$, \  $\alpha_j
\in\N$ \quad $P_{\overrightarrow\alpha}^n (x) = \displaystyle{1\over
\vmid{\overrightarrow\alpha}!} \  D^{\overrightarrow\alpha} \
x^{n+\vmid{\overrightarrow\alpha}}$

then \ $P_{\overrightarrow\alpha}^n$ \ is of degre $n$
$$\eqalign{
\vmid{\alpha} ! \ P_{\overrightarrow\alpha}^{\vmid{\alpha}-1} (x) &= {(2 \
\vmid{\alpha} - \alpha_0 - 1)!\over (2\vmid{\alpha} - 1)!} P_\alpha (x) \cr
{\partial\over \partial x_0} \ P_{\overrightarrow\alpha}^n (x) &= n \
P_{\overrightarrow\alpha}^{n-1} (x)\cr
{\partial\over \partial x_k} \ P_{\overrightarrow\alpha}^n (x) &=
P_{(\alpha_{1},\ldots , \alpha_{k}+1,\ldots , \alpha_{2m+1})}^{n-1} (x). \cr}$$

\vskip 1,5cm
{\bf 6. The Cauchy kernel of holomorphic Cliffordian functions}
\bigskip\medskip
\hskip0,45cm Following Brackx, Delanghe and Sommen [1], recall that there
exists a Cauchy kernel connected with the theory of monogenic functions. In our
situation, when we study functions of the type~:
$$f : S \oplus V \fle \R_{0,2m+1},$$
the related Cauchy kernel is :
$$E(x) = {1\over \omega_m} \ {x^*\over \vmid{x}^{2m+2}}, \quad x\in S \oplus V
\setminus \{ 0\} ,\leqno (7)$$
where \ $\omega_m = 2 \pi^{m+1} \ \displaystyle{1\over \Gamma (m+1)}$ \ is the
area of the unit sphere in \ $\R^{2m+2}$.

\medskip
\hskip0,45cm Recall also that \ $E(x)$ \ is a monogenic function with singularity
at the origin, i.e~:
$$DE(x) = \delta \quad\hbox{for}\quad x \in S \oplus V$$

where $\delta$ is the Dirac measure.

Let \ $\omega (y) = dy_0 \wedge \cdots \wedge dy_{2m+1}$ \ and \ $\gamma (y) =
\displaystyle\sum_{i=0}^{2m+1} (-1)^i \ e_i dy_0 \wedge \cdots \wedge
\widehat{dy_i} \wedge \cdots \wedge dy_{2m+1}$.

\bigskip\medskip
\hskip0,45cm Then, we have :

\bigskip
{\bf Theorem.-} \ [Integral representation formula (general case)] \  [1].

If \ $f\in {\cal C}^1 (U, \ \R_{0,2m+1})$,  then~:

$$\int_{\partial\Omega} E(y-x)  \gamma (y) f(y) - \int_\Omega E(y-x) Df(y) \omega
(y) = \cases{f(x), &$x\in \buildrel\, \, \circ\over \Omega$\cr
0, &$x\notin \Omega$,\cr}$$

where $\Omega$ is an oriented compact differentiable variety of dimension $2m+2$
with boundary $\partial \Omega$ \ and \ $\Omega \subset U$.

\medskip
\hskip0,45cm From this theorem follows the following integral representation
formula for monogenic functions called also the Cauchy representation formula
[1]~:

\bigskip\bigskip
{\bf Theorem.-} \  If $f$ is monogenic in $U$ and if \ $\Omega\subset U$
$$\int_{\partial \Omega} E(y-x)\gamma (y) f(y)  = \cases{f(x),
&$x\in\buildrel\,\,
\circ\over \Omega$\cr
0, &$x\notin \Omega$\cr}$$

\hskip0,45cm It is natural to have an integral representation formula of this type
because the Cauchy operator $D$ is of order 1. In our situation, the operator \
$D\Delta^m$ \ which gives the holomorphic Cliffordian functions is of order $2m+1$
and the corresponding integral formula would be  much more complicated.

\hskip0,45cm But, the first step to obtain such a formula, is to exhibit an
analogous of the Cauchy kernel.

\hskip0,45cm Remember that the fundamental solution of the iterated Laplacian,
i.e. the function \ $h : S \oplus V \setminus \{ 0\} \fle \R$  verifying the
equation \ $\Delta^{m+1}h(x) = 0$ \ for \ $x\in S \oplus V \setminus \{ 0\}$,  is
in fact well-known : that is
$$h(x) = \ell n \ \vmid{x}, \quad x\in S \oplus V \setminus \{ 0\} .$$

\hskip0,45cm Recall briefly the idea : using spherical coordinates, i.e.
introducing  \ $\rho = \ \vmid{x}$,  the radial form of the Laplacian is
$$\Delta_\rho = {d^2\over d\rho^2} + {2m+1\over \rho} \ {d\over d\rho} .$$

Calculating the iterated Laplacian, one get, for $k\in\N$~:
$$\Delta^k_\rho \ell n \ \rho = (-1)^{k+1} 2^{k-1} (k-1)! (2m) (2m-2) \cdots (2m 
- 2k+2) \ {1\over \rho^{2k}}$$
and, thus, when \ $k = m+1$, \ one has outside the singularity :
$$\Delta_\rho^{m+1} \ell n \ \rho = 0.$$

\medskip
\hskip0,45cm Similarly as in the complex case when we know that \ $\ell n
\sqrt{x^2+y^2}$ \ is the fundamental solution of the Laplace equation and when
we write it as \ $\displaystyle{1\over 2} \ \ell n \ (z\overline z)$,  here
also we will resort to the relation \ $xx^* = \ \vmid{x}^2$ \ for \ $x = x_0 +
\vec x \in S \oplus V $ \ and the final conclusion of our first step
is~:

\hskip0,45cm The fundamental solution of the iterated Laplacian \ $\Delta^{m+1}$ \
is \ $h(x) = \displaystyle{1\over 2} \ \ell n \ (xx^*)$.

\medskip
\hskip0,45cm Now, according to the property (ii) of \SS 3. \ $h(x)$ being a
polyharmonic function of order \ $m+1$, then \ $D^* \bigl({1\over 2} \ell n \
(xx^*)\bigr)$ \ will be a holomorphic Cliffordian function on \ $S \oplus V
\setminus \{ 0\}$.  But
$$D^* \ \bigl( {1\over 2} \ell n \ (xx^*)\bigr) = {1\over 2} \
{D^*(\vmid{x}^2)\over \vmid{x}^2} = {x^*\over \vmid{x}^2} = x^{-1}.$$

\medskip
\hskip0,45cm In this way, we have found the first holomorphic Cliffordian function
with singularity at the origin.

\medskip
\hskip0,45cm Again, according to the remark of \SS 2. since \ $x^{-1}$ \ is
holomorphic Cliffordian on \ $S \oplus V \setminus \{ 0\}$,  then \
$\Delta^m(x^{-1})$ \ should be monogenic on the same set.  Let us compute
$$\Delta^m (x^{-1}) = \Delta^m D^* \ell n \ \rho = D^* \Delta^m \ \ell n \ \rho ,$$
where we have noted \ $\rho = \ \vmid{x} = (xx^*)^{{1\over 2}}$.

Now explicitly,
$$\eqalign{ \Delta^m (x^{-1})
&= D^* (-1)^{m+1} 2^{m-1} (m-1)! (2m) (2m-2) \cdots 2 . {1\over \rho^{2m}} =\cr
&= (-1)^{m+1} 2^{2m-1} (m-1)! m! \ D^* \ \Bigl( {1\over \rho^{2m}}\Bigr) = \cr
&= (-1)^m 2^{2m-1} (m!)^2 \ {1\over {(\vmid{x}^2)}^{m+1}} \ D^* (\vmid{x}^2) =\cr
&= (-1)^m 2^{2m} (m!)^2 \ {x^*\over \vmid{x}^{2m+2}} =\cr
&= (-1)^m 2^{2m} (m!)^2 \ \omega_m \  E(x).\cr}$$

Thus, we get :
$${(-1)^m (m+1)\over 2^{2m+1} m! \ \pi^{m+1}} \ \Delta^m (x^{-1}) = E(x).$$

It becomes natural to introduce a new kernel :
$$N(x) = \varepsilon_m \ x^{-1},$$

where \ $\varepsilon_m = (-1)^m \ \displaystyle{m+1\over 2^{2m+1} m! \ \pi^{m+1}}$.

\medskip
Remember the basic properties of the kernel \ $N(x)$~:

\medskip
\hskip0,45cm (i) \quad It is related to the Cauchy kernel of the monogenic
functions $E$ by~:
$$\Delta^m N(x) = E(x), \quad x\in S \oplus V \setminus \{ 0\}.$$

\bigskip
\hskip0,45cm (ii) \ \ $N$ is holomorphic Cliffordian on \ $S\oplus V \setminus \{
0\}$  because~:
$$D \Delta^m N(x) = DE(x) = \delta .$$

\vskip 1,5cm
{\bf 7. Integral representation formula for holomorphic Cliffordian functions}
\bigskip
\hskip0,45cm Let \ $f : S \oplus V \fle R_{0,2m+1}$ \ be a function of class \
${\cal C}^{2m+1}$ \ and $B$ be the unit ball in \ $\R^{2m+2}$.  According to
[1], for \  $x\in\buildrel\, \, \circ\over B$, we have
$$f(x) = \int_{\partial B} E(y-x) \gamma (y) f(y)  - \int_B E(y-x)
Df(y) \omega (y).$$

Substitute \ $\Delta^mN$ \ on the place of $E$, one has~:
$$f(x) = \int_{\partial B} \Delta^m N(y-x) \gamma (y) f(y)  - \int_B
\Delta^m N(y-x) Df(y) \omega (y).$$

Making use of the Green's formula :
$$\int_\Omega u\Delta v = \int_\Omega v\Delta u + \int_{\partial\Omega} u \
{\partial v\over \partial n} - \int_{\partial \Omega} v \ {\partial u\over \partial
n}$$

applied on the second integral with \ $u = Df$ \ and \ $v = \Delta^{m-1} N$, \ we
will deduce~:

$$\eqalign{
&f(x) = \int_{\partial B} \Delta^m N(y-x) \gamma (y) f(y)  - \int_B
\Delta^{m-1} N(y-x) \Delta Df(y) \omega (y)\cr
&- \int_{\partial B} \ \Bigl( {\partial\over \partial n} \Delta^{m-1}
N(y-x)\Bigr) Df(y) d\sigma_y + \int_{\partial B} \ \bigl( \Delta^{m-1}
N(y-x)\bigr) \ {\partial\over \partial n} Df(y) d\sigma_y .\cr}$$

\medskip
\hskip0,45cm Iterating the process of applying the Green's formula on the second
integral of the preceding formula with \ $u = D\Delta f$ \ and \ $v =
\Delta^{m-2}N$,  we will deduce a sum of six integrals as follows~:
$$\eqalign{f(x) &= 
\int_{\partial B} (\Delta^m N) \gamma f - \int_B (\Delta^{m-2}N) D\Delta^2f -\cr
&-\int_{\partial B} \ \Bigl( {\partial\over\partial n} \Delta^{m-2}
N\Bigr) D\Delta f + \int_{\partial B} \ (\Delta^{m-2} N) \ {\partial\over \partial
n} D\Delta f\cr
&- \int_{\partial B} \ \Bigl( {\partial\over \partial n} \Delta^{m-1} B\Bigr) Df +
\int_{\partial B} (\Delta^{m-1}N) \ {\partial\over \partial n} Df.\cr}$$

So, applying \ $m$ \ times the Green's formula, we have :
$$\eqalign{
f(x)
&= \int_{\partial B} \bigl(\Delta^m N(y-x)\bigr) \gamma (y) f(y) \cr
&- \sum_{k=1}^m \ \int_{\partial B} \ \Bigl( {\partial\over \partial n}
\Delta^{m-k} N(y-x)\Bigr) D\Delta^{k-1} f(y) d\sigma_y\cr
&+ \sum_{k=1}^m \ \int_{\partial B} \ \bigl( \Delta^{m-k} N(y)x)\bigr) \
{\partial\over \partial n} D\Delta^{k-1} f(y) d\sigma_y\cr
&- \int_B N(y-x) D\Delta^m f(y) \omega (y).\cr}$$

This would be the general integral representation formula for functions \ $f : S
\oplus V \fle \R_{0,2m+1}$.

\medskip
\hskip0,45cm The Cauchy integral formula for holomorphic Cliffordian functions
will be obtained erasing the last integral because in that case \ $D\Delta^mf = 0$.

\medskip
\hskip0,45cm Remark that the obtained Cauchy integral formula involves \ $2m+1$ \
integrals on \ $\partial B$.  That means that, for holomorphic Cliffordian
function, one can reconstitute the values of \ $f$ \ in a point of the interior of
\ $B$ \ knowing the values on \ $\partial B$ \ of $f$, \ $D\Delta^{k-1}f$ \ and \
$\displaystyle{\partial\over\partial n} D\Delta^{k-1} f$, \ with \ $k = 1,\ldots
,m$.

\medskip
\hskip0,45cm Remark also, that when \ $m = 0$, \ i.e.  the case of holomorphic
functions, we have~:
$$N(z) = {1\over 2\pi} \ \cdot \ {1\over z}, \quad E(z) = {1\over 2\pi} \
{\overline z\over \vmid{z}^2} .$$

\vskip 1,5cm
{\bf 8. Taylor expansion of a holomorphic Cliffordian function}
\bigskip
\hskip0,45cm Here we will imitate the well-know process for the obtention of a
Taylor formula for holomorphic functions starting with the Cauchy formula and
developping the Cauchy kernel. Our Cauchy kernel is~:
$$N(y-x) = \varepsilon_m (y-x)^{-1}.$$

In order to developp \ $(y-x)^{-1}$, \ let us proceed as follows~:
$$\displaylines{
(y-x)^{-1} = {\Bigl( y (1-y^{-1} x)\Bigr)}^{-1} = (1-y^{-1}x)^{-1} y^{-1} =\cr
= y^{-1} + y^{-1} xy^{-1} + y^{-1}xy^{-1}xy^{-1} + \cdots + {(y^{-1}x)}^n y^{-1} +
\cdots\cr}$$

In view of \ $yy^* = \  \vmid{y}^2$, \ we have \ $y^{-1} = \displaystyle{y^*\over
\vmid{y}^2}$,  and thus~:
$$(y-x)^{-1} = \sum_{n=0}^\infty \ {{(y^*x)}^n y^*\over \vmid{y}^{2n+2}} .$$

Let have a look at the second term of this developpement~:
$$\eqalign{ 
&y^*xy^* = (y_0 - \vec y\ ) x (y_0 - \vec y\ ) =\cr
&= (e_0xe_0) y_0^2 + \sum_{j=1}^{2m+1} (e_0x e_j) y_0 (-y_j) + \sum_{k=1}^{2m+1}
(e_kxe_0) (-y_k)y_0 + \sum_{j,k=1}^{2m+1} \ (e_jxe_k) y_jy_k.\cr}$$

It is not difficult to observe that the polynomials \ $P_\alpha (x)$ \ appear
again and one can writte~:
$$y^* xy^* = \sum_{\vmid{\alpha} = 2} \ P_\alpha (x) Y^\alpha ,$$

where we have made use of the notation :
$$Y^\alpha = y_0^{\alpha_{0}} (-y_1)^{\alpha_{1}} \cdots
(-y_{2m+1})^{\alpha_{2m+1}}.$$

A straightforward calculation gives finally :
$$(y-x)^{-1} = \sum_{k=1}^\infty \ {1\over \vmid{y}^{2k}} \ \sum_{\vmid{\alpha}
= k} \ P_\alpha (x) Y^\alpha$$
or more concisely :
$$(y-x)^{-1} = \sum_{\vmid{\alpha}=1}^\infty \ P_\alpha (x) \ {Y^\alpha\over
\vmid{y}^{2\vmid{\alpha}}} .$$

In order to obtain the Taylor series of $f$, take the Cauchy integral formula and
substitute the expansion of $N(y-x)$.

Observe that \ $\Delta_x^m N(y-x) = \Delta_y^m  N(y-x)$,  so that in the first
integral of the Cauchy formula, we have~:
$$\eqalign{
&\int_{\partial B} \Delta^m N(y-x) \gamma (y) f(y)  = \int_{\partial B}
\Delta_y^m
\ \Bigl( \varepsilon_m \ \sum_{\vmid{\alpha} = 1}^\infty \ P_\alpha (x) \
{Y^\alpha\over \vmid{y}^{2\vmid{\alpha}}}\Bigr) \gamma (y) f(y)  =\cr
&=\varepsilon_m \ \sum_{\vmid{\alpha}=1}^\infty P_\alpha(x) \int_{\partial B} \
\Bigl( \Delta_y^m \ {Y^\alpha\over \vmid{y}^{2\vmid{\alpha}}}\Bigr) \gamma (y)
f(y)
 = \sum_{\vmid{\alpha} = 1}^\infty P_\alpha (x) \ A_\alpha^{(0)},\cr}$$

where the \ $A_\alpha^{(0)}$  are in \ $\R_{0,2m+1}$ \ and are given by~:

$$A_\alpha^{(0)} = \varepsilon_m \int_{\partial B} \ \Bigl( \Delta_y^m
\ {Y^\alpha\over \vmid{y}^{2\vmid{\alpha}}}\Bigr) \gamma (y)\  f(y).$$

Similarly, as the other integrals in the Cauchy formula are concerned, we have~:

$${\partial\over \partial n} \ \Delta^\ell N(y-x) = {\partial\over \partial n_y} \
\Delta_y^\ell N(y-x)$$

which allows to deduce finally : 
$$f(x) = \sum_{\vmid{\alpha} = 1}^\infty P_\alpha (x) C_\alpha ,$$

where the coefficients \ $C_\alpha \in \R_{0,2m+1}$,  and more precisely~: 
$$C_\alpha = A_\alpha^{(0)} + A_\alpha^{(1)} + \cdots + A_\alpha^{(2m)}$$

with :
$$A_\alpha^{(j)} = \varepsilon_m \int_{\partial S} \ \Bigl( {\partial\over
\partial n_y} \Delta_y^{m-j} \ {Y^\alpha\over \vmid{y}^{2\vmid{\alpha}}}\Bigr)
D\Delta^{j-1} f(y) d\sigma_y, \ j = 1,\ldots ,m$$

and
$$A_\alpha^{(\ell + m)} = \varepsilon_m \int_{\partial S} \ \Bigl(
\Delta_y^{m-\ell} \ {Y^\alpha\over \vmid{y}^{2\vmid{\alpha}}}\Bigr) \
{\partial\over \partial n} D\Delta^{\ell -1} f(y) d\sigma_y, \ \ell = 1,\ldots ,m.$$

At the end of this paragraph let us prove that the polynomials \ $P_\alpha$ \ span
the space of  polynomial solutions of \ $D\Delta^m f = 0$.  Indeed, according to
the Taylor expansion if \ $P(x)$ \ is an arbitrary polynomial, we have~:
$$P(x) = \sum_{\vmid{\alpha} = 1}^\infty \ P_\alpha (x) C_\alpha$$

as a holomorphic Cliffordian function. But \ $P$ \ is a polynomial, so that the
sum is finite~:
$$P(x) = \sum_{\vmid{\alpha} = 1}^d \ P_\alpha (x) \ C_\alpha$$

and this shows that $P$ is a linear  (right) combination of the $P_\alpha$.

\bigskip
Let $Q$ be any polynomial of degree smaller or equal to $2m$,  then $Q$ is
holomorphic Cliffordian and
$$Q(x) = \sum_{\vmid{\alpha} = 1}^{2m+1} P_\alpha (x) C_\alpha.$$

\vskip 1,5cm
{\bf 9. Laurent series}
\bigskip
\hskip0,45cm Consider a function which is holomorphic Cliffordian on a punctured
neighborhood of the origin, say, for example on \ $B \setminus \{ 0\}$,  where $B$
is the unit ball in \ $S \oplus V$.

\hskip0,45cm Let \ $\Gamma_1$ \ and \ $\Gamma_2$ \ be two balls, centered at the
origin, with radii  \ $r_1$ and $r_2$, respectively, and such that \ $0 < r_1 <
r_2 < 1$.  One can applied the Cauchy  representation formula on the region, which
is limited by $\Gamma_1$ and $\Gamma_2$,  namely on \ $\Gamma_2 \setminus
\Gamma_1$. Those integrals, taken on \ $\partial\Gamma_2$,  will give us, as in
the previous paragraph, the regular part of the Laurent series. Because of the
sense of the integration, we have now to integrate on \ $\partial\Gamma_1$ \
those terms of the representation formula, which contain \ $N(x-y)$ \ and its
derivatives.

\medskip
\hskip0,45cm In this way, one needs to developp \ $(x-y)^{-1}$. So~:
$$\eqalign{ (x-y)^{-1}
&= {(x (1-x^{-1}y))}^{-1} = (1-x^{-1}y)^{-1} x^{-1} =\cr
&= x^{-1} + x^{-1}y x^{-1} + x^{-1}yx^{-1}yx^{-1} + \cdots \cr
&\ \cdots + (x^{-1} y)^k \ x^{-1} + \cdots \cr
&= x^{-1} + \ \sum_{i=0}^{2m+1} \ (x^{-1} e_i x^{-1}) y_i  \ + \cr
&+ \sum_{0\leq i_{1}, i_{2} \leq 2m+1} \ (x^{-1} e_{i_{1}} x^{-1} e_{i_{2}}
x^{-1} + x^{-1} e_{i_{2}} x^{-1} e_{i_{1}} x^{-1}) y_{i_{1}} y_{i_{2}} + \cdots
\cr}$$

\medskip
\hskip0,45cm Remark that the rational functions appearing in the last
developpement are of negative powers on $x$, \ $\resp.$~: -1, -2, -3,$\ldots$ \ .

\medskip
\hskip0,45cm Using a similar manner of notation as in the case of the polynomials
\ $P_\alpha (x)$,  we set \ $\beta = (\beta_0, \beta_1,\ldots , \beta_{2m+1})$, 
with \ $\beta_i\in\N$ \ and \ $\vmid{\beta} = \displaystyle\sum_{i=0}^{2m+1}
\beta_i$.  Consider again the set \ $\{ e_\nu \}$,  where \ $e_0$ \ is written \
$\beta_0$ \ times, \ $e_1, \beta_1$ times, etc $\ldots$ \ and \ $e_{2m+1},
\beta_{2m+1}$ \ times. Set now~:
$$S_\beta (x) = \sum_{\hbox{\got S}} \ \prod_{\nu = 1}^{\vmid{\beta}} \ (x^{-1}
e_{\sigma (\nu )}) x^{-1},$$

the sum being expanded over all distinguishable elements $\sigma$ of the permutation
group  \got S .

\hskip0,45cm $S_\beta (x)$ \ is left and right holomorphic Cliffordian.

\medskip
\hskip0,45cm We recognize easily \ $S_{(1,0,0,0)} (x) = x^{-1} e_0 x^{-1}$, \
$S_{(0,1,0,0)} (x) = x^{-1} e_1 x^{-1}$, \break $S_{(0,1,1,0)} (x) = x^{-1} e_1
x^{-1} e_2 x^{-1} + x^{-1} e_2 x^{-1} e_1 x^{-1}$ \ in the special case when $m=1$.
Remark also that \ $S_0(x) = x^{-1}$ \ and that the power of \ $x^{-1}$ \ in \
$S_\beta$ \ is exactly \ $\vmid{\beta} + 1$. 

\medskip
\hskip0,45cm Thus, we get :
$$N(x-y) = \varepsilon_m \ \sum_{\vmid{\beta} = 0}^\infty \ S_\beta (x) y^\beta ,$$

where \ $y^\beta = (y_0)^{\beta_{0}} (y_1)^{\beta_{1}} \ldots
(y_{2m+1})^{\beta_{2m+1}}$.

\bigskip
\hskip0,45cm In the same way as in paragraph 8, one deduces the following Laurent
series for a holomorphic Cliffordian function \ $f : B \setminus \{ 0\} \fle
\R_{0,2m+1}$, \ $B \subset S \oplus V$~: for each \ $x\in B \setminus \{ 0\}$,  we
have
$$f(x) = \sum_{\vmid{\beta}=0}^\infty \ S_\beta(x) D_\beta + \sum_{\vmid{\alpha} =
1}^\infty P_\alpha (x) C_\alpha ,$$

where \ $C_\alpha$ \ and \ $D_\beta$ \ belong to \ $\R_{0,2m+1}$.

\medskip
\hskip0,45cm The first sum is the analogous of the singular part of a Laurent
expansion for a holomorphic function, while the second sum represents the
analogous of its regular part.

\medskip
\hskip0,45cm Here, we centered our expansions at the origin. Of course, they remain
valid in neighborhoods of every point \ $a\in S \oplus V$. If \ $f : B \setminus
\{ a\} \fle \R_{0,2m+1}$, \ is a holomorphic Cliffordian function, where $B$ is a
ball centered at $a$,  then for every $x\in B \setminus \{ a\}$,  one has~:
$$f(x) = \sum_{\vmid{\beta} = 0}^\infty \ S_\beta (x-a) D_\beta +
\sum_{\vmid{\alpha} = 1}^\infty \ P_\alpha (x-a) C_\alpha ,$$

with \ $C_\alpha , D_\beta \in \R_{0, 2m+1}$.

\bigskip\bigskip
{\bf Remark :} \ the rational functions \ $S_\beta (x)$ \ are left and right
holomorphic Cliffordian.

\bigskip\medskip
\hskip 0,45cm The present paper is a detailed exposition of part of the results
announced in [8]. However, some modifications were brought, especially
concerning the multiplicative constants appearing in the definitions of the
polynomials \
$P_\alpha (x)$ \ and the rational functions $S_\beta (x)$.

\vfill\eject
\centerline{\bf Bibliographie}
\bigskip
\bb 1&F. BRACKX, R. DELANGHE, F. SOMMEN&Clifford analysis,&Pitman, (1982) & &

\bb 2&C.A. DEAVORS&The quaternion calculus&Am. Math. Monthly. (1973), 995-1008& &

\bb 3&R. DELANGHE, F. SOMMEN, V. SOUC\v EK&Clifford Algebra and
Spinor-valued functions&Kluwer Academic Publishers& &

\bb4&R. FUETER&Die Funktionnentheorie der Differentialgleichungen \ $\Delta
u = 0$ und $\Delta\Delta u = 0$ mit vier reellen Variablen.& Comment Math. Helv 7
(1935), 307-330& &

\bb5&R. FUETER&Uber die analytische Darstellung der regularen
Funktionen einer Quaternionenvariabelen&Comm. Math. Helv.8 (1936), 371-378& &

\bb6&G. LAVILLE&Une famille de solutions de l'\'equation de Dirac avec champ
\'electromagn\'etique quelconque&C.R. Acad. Sci. Paris t. 296 (1983), 1029-1032& &

\bb7&G. LAVILLE&Sur l'\'equation de Dirac avec champ \'electromagn\'etique
quelconque&Lectures Notes in Math. 1165, Springer-Verlag (1985), 130-149& &

\bb8&G. LAVILLE, I. RAMADANOFF&Fonctions holomorphes Cliffordiennes&C.R. Acad,
Sc. Paris, 326, s\'erie I (1998), 307-310& &

\bb9&H. MALONEK&Powers series representation for monogenic functions in \
$\R^{n+1}$ based on a permutational product&Complex variables, vol 15 (1990),
181-191& &

\bb10&V.P. PALAMODOV&On ``holomorphic" functions of several quaternionic \break
variables&C.A. Aytama (ed) Linear topological spaces and complex analysis II,
Ankara (1995), 67-77& &

\bb11&L. PERNAS&Holomorphic quaternionienne&preprint, (1997)& &

\vskip1cm

{\parindent=8cm
\item{UPRES-A 6081} D\'epartement de Math\'ematiques
\item{}Universit\'e de Caen
\item{}14032 CAEN Cedex France

\medskip
\item{}glaville@math.unicaen.fr
\item{}rama@math.unicaen.fr
\par}
}

\end